\magnification=\magstep1


\def\item{\vskip1.3pt\hang\textindent}


\tolerance=300
\pretolerance=200
\hfuzz=1pt
\vfuzz=1pt

\hoffset 0cm            
\hsize=5.8 true in
\vsize=9.5 true in

\def\rightheadline{\hfil\smc\lastname\hfil\tenbf\folio}
\def\leftheadline{\tenbf\folio\hfil\smc\lastname\hfil}
\headline={\ifodd\pageno\rightheadline\else\leftheadline\fi}
\newdimen\dimenone
\def\checkleftspace#1#2#3#4#5{
 \dimenone=\pagetotal
 \advance\dimenone by -\pageshrink   
 \ifdim\dimenone>\pagegoal          
   \else\dimenone=\pagetotal
        \advance\dimenone by \pagestretch
        \ifdim\dimenone<\pagegoal
          \dimenone=\pagetotal
          \advance\dimenone by#1         
          \setbox0=\vbox{#2\parskip=0pt                
                       \hyphenpenalty=10000
                       \rightskip=0pt plus 5em
                       \noindent#3 \vskip#4}    
        \advance\dimenone by\ht0
        \advance\dimenone by 3\baselineskip
        \ifdim\dimenone>\pagegoal\vfill\eject\fi
          \else\eject\fi\fi}

\parindent=35pt
\mathsurround=1pt
\parskip=1pt plus .25pt minus .25pt
\normallineskiplimit=.99pt

\mathchardef\emptyset="001F 

\def\Int{\mathop{\rm int}\nolimits}
%



\def\1{{\bf1}}\def\0{{\bf0}}

\def\({\bigl(}  \def\){\bigr)}
\def\<{\mathopen{\langle}}\def\>{\mathclose{\rangle}}

\def\Z{{\mathchoice{{\hbox{$\rm Z\hskip 0.26em\llap{\rm Z}$}}}%
{{\hbox{$\rm Z\hskip 0.26em\llap{\rm Z}$}}}%
{{\hbox{$\scriptstyle\rm Z\hskip 0.31em\llap{$\scriptstyle\rm Z$}$}}}{{%
\hbox{$\scriptscriptstyle\rm Z$\hskip0.18em\llap{$\scriptscriptstyle\rm Z$}}}}}}

\def\N{{\mathchoice{\hbox{$\rm I\hskip-0.14em N$}}%
{\hbox{$\rm I\hskip-0.14em N$}}%
{\hbox{$\scriptstyle\rm I\hskip-0.14em N$}}%
{\hbox{$\scriptscriptstyle\rm I\hskip-0.10em N$}}}}

\def\R{{\mathchoice{\hbox{$\rm I\hskip-0.14em R$}}%
{\hbox{$\rm I\hskip-0.14em R$}}%
{\hbox{$\scriptstyle\rm I\hskip-0.14em R$}}%
{\hbox{$\scriptscriptstyle\rm I\hskip-0.10em R$}}}}

\def\K{{\mathchoice{\hbox{$\rm I\hskip-0.15em K$}}%
{\hbox{$\rm I\hskip-0.15em K$}}%
{\hbox{$\scriptstyle\rm I\hskip-0.15em K$}}%
{\hbox{$\scriptscriptstyle\rm I\hskip-0.11em K$}}}}

\def\qed{\hfill {\hbox{[\hskip-0.05em ]}}}

\def\C{{\mathchoice%
{\hbox{$\rm C\hskip-0.47em\hbox{%
\vrule height 0.58em width 0.06em depth-0.035em}$}\;}%
{\hbox{$\rm C\hskip-0.47em\hbox{%
\vrule height 0.58em width 0.06em depth-0.035em}$}\;}%
{\hbox{$\scriptstyle\rm C\hskip-0.46em\hbox{%
$\scriptstyle\vrule height 0.365em width 0.05em depth-0.025em$}$}\;}
{\hbox{$\scriptscriptstyle\rm C\hskip-0.41em\hbox{
$\scriptscriptstyle\vrule height 0.285em width 0.04em depth-0.018em$}$}\;}}}

\def\.{{\cdot}}
\def\|{\Vert}
\def\ssk{\smallskip}
\def\msk{\medskip}
\def\bsk{\bigskip}
\def\giantskip{\vskip2\bigskipamount}

\def\giantbreak{\par \ifdim\lastskip<2\bigskipamount \removelastskip
         \penalty-400 \giantskip\fi}

\def\nin{\noindent}
\def\cen{\centerline}
\def\pagebreak{\vskip 0pt plus 0.0001fil\break}
\def\linebreak{\break}

\def\epsilon{\varepsilon}

\font\ninerm=cmr9
\font\eightrm=cmr8
\font\sixrm=cmr6

\font\eightbf=cmbx8
\font\sixbf=cmbx6

\font\eighti=cmmi8
\font\sixi=cmmi6
\font\ninesy=cmsy9
\font\eightsy=cmsy8
\font\sixsy=cmsy6

\font\eightit=cmti8


\font\eightsl=cmsl8

\font\eighttt=cmtt8
\font\bfone=cmbx10 scaled\magstep1 
\font\smc=cmcsc10

\font\small=cmcsc8

\def\no #1. {\bigbreak\vskip-\parskip\noindent\bf #1. \quad\rm}

\def\Proposition #1. {\checkleftspace{0pt}{\bf}{Theorem}{0pt}{}
\bigbreak\vskip-\parskip\noindent{\bf Proposition #1.}
\quad\it}

\def\Theorem #1. {\checkleftspace{0pt}{\bf}{Theorem}{0pt}{}
\bigbreak\vskip-\parskip\noindent{\bf  Theorem #1.}
\quad\it}
\def\Corollary #1. {\checkleftspace{0pt}{\bf}{Theorem}{0pt}{}
\bigbreak\vskip-\parskip\nin{\bf Corollary #1.}
\quad\it}
\def\Lemma #1. {\checkleftspace{0pt}{\bf}{Theorem}{0pt}{}
\bigbreak\vskip-\parskip\noindent{\bf  Lemma #1.}\quad\it}

\def\Definition #1. {\checkleftspace{0pt}{\bf}{Theorem}{0pt}{}
\rm\bigbreak\vskip-\parskip\noindent{\bf Definition #1.}
\quad}

\def\Remark #1. {\checkleftspace{0pt}{\bf}{Theorem}{0pt}{}
\rm\bigbreak\vskip-\parskip\noindent{\bf Remark #1.}\quad}

\def\Exercise #1. {\checkleftspace{0pt}{\bf}{Theorem}{0pt}{}
\rm\bigbreak\vskip-\parskip\noindent{\bf Exercise #1.}
\quad}

\def\Example #1. {\checkleftspace{0pt}{\bf}{Theorem}{0pt}{}
\rm\bigbreak\vskip-\parskip\noindent{\bf Example #1.}\quad}
\def\Examples #1. {\checkleftspace{0pt}{\bf}{Theorem}{0pt}
\rm\bigbreak\vskip-\parskip\noindent{\bf Examples #1.}\quad}

\newcount\problemnumb \problemnumb=0
\def\Problem{\global\advance\problemnumb by 1\bigbreak\vskip-\parskip\noindent
{\bf Problem \the\problemnumb.}\quad\rm }

\def\Proof#1.{\rm\par\ifdim\lastskip<\bigskipamount\removelastskip\fi\smallskip
            \noindent {\bf Proof.}\quad}

\nopagenumbers

\def\author{}
\def\lastname{}
\def\thanks#1{\footnote*{\eightrm#1}}
\def\title{}

\def\lastname{}
\def\h{{\textstyle{1\over2}}}

\def\ep{\epsilon}

\def\text{\textstyle} 
\def\disp{\displaystyle} 
\def\d{{\,\rm d}}

\def\and{{\rm and }}

\expandafter\edef\csname amssym.def\endcsname{%
       \catcode`\noexpand\@=\the\catcode`\@\space}
\catcode`\@=11
\def\undefine#1{\let#1\undefined}
\def\newsymbol#1#2#3#4#5{\let\next@\relax
 \ifnum#2=\@ne\let\next@\msafam@\else
 \ifnum#2=\tw@\let\next@\msbfam@\fi\fi
 \mathchardef#1="#3\next@#4#5}
\def\mathhexbox@#1#2#3{\relax
 \ifmmode\mathpalette{}{\m@th\mathchar"#1#2#3}%
 \else\leavevmode\hbox{$\m@th\mathchar"#1#2#3$}\fi}
\def\hexnumber@#1{\ifcase#1 0\or 1\or 2\or 3\or 4\or 5\or 6\or 7\or 8\or
 9\or A\or B\or C\or D\or E\or F\fi}

\font\tenmsb=msbm10
\font\sevenmsb=msbm7
\font\fivemsb=msbm5
\newfam\msbfam
\textfont\msbfam=\tenmsb
\scriptfont\msbfam=\sevenmsb
\scriptscriptfont\msbfam=\fivemsb
\edef\msbfam@{\hexnumber@\msbfam}
\def\Bbb#1{{\fam\msbfam\relax#1}}

\newsymbol\Bbbk 207C
\def\widehat#1{\setbox\z@\hbox{$\m@th#1$}%
 \ifdim\wd\z@>\tw@ em\mathaccent"0\msbfam@5B{#1}%
 \else\mathaccent"0362{#1}\fi}
\def\widetilde#1{\setbox\z@\hbox{$\m@th#1$}%
 \ifdim\wd\z@>\tw@ em\mathaccent"0\msbfam@5D{#1}%
 \else\mathaccent"0365{#1}\fi}
\font\teneufm=eufm10
\font\seveneufm=eufm7
\font\fiveeufm=eufm5
\newfam\eufmfam
\textfont\eufmfam=\teneufm
\scriptfont\eufmfam=\seveneufm
\scriptscriptfont\eufmfam=\fiveeufm

\catcode`@=11 

\expandafter\edef\csname amssym.def\endcsname{%
       \catcode`\noexpand\@=\the\catcode`\@\space}
\font\eightmsb=msbm8
\font\sixmsb=msbm6
\font\fivemsb=msbm5
\font\eighteufm=eufm8
\font\sixeufm=eufm6
\font\fiveeufm=eufm5
\newskip\ttglue
\def\eightpoint{\def\rm{\fam0\eightrm}%
  \textfont0=\eightrm \scriptfont0=\sixrm \scriptscriptfont0=\fiverm
  \textfont1=\eighti \scriptfont1=\sixi \scriptscriptfont1=\fivei
  \textfont2=\eightsy \scriptfont2=\sixsy \scriptscriptfont2=\fivesy
  \textfont3=\tenex \scriptfont3=\tenex \scriptscriptfont3=\tenex
\textfont\eufmfam=\eighteufm
\scriptfont\eufmfam=\sixeufm
\scriptscriptfont\eufmfam=\fiveeufm
\textfont\msbfam=\eightmsb
\scriptfont\msbfam=\sixmsb
\scriptscriptfont\msbfam=\fivemsb
  \def\it{\fam\itfam\eightit}%
  \textfont\itfam=\eightit
  \def\sl{\fam\slfam\eightsl}%
  \textfont\slfam=\eightsl
  \def\bf{\fam\bffam\eightbf}%
  \textfont\bffam=\eightbf \scriptfont\bffam=\sixbf
   \scriptscriptfont\bffam=\fivebf
  \def\tt{\fam\ttfam\eighttt}%
  \textfont\ttfam=\eighttt
  \tt \ttglue=.5em plus.25em minus.15em
  \normalbaselineskip=9pt
  \def\MF{{\manual opqr}\-{\manual stuq}}%
  \let\big=\eightbig
  \setbox\strutbox=\hbox{\vrule height7pt depth2pt width\z@}%
  \normalbaselines\rm}
\def\eightbig#1{{\hbox{$\textfont0=\ninerm\textfont2=\ninesy
  \left#1\vbox to6.5pt{}\right.\n@space$}}}


\csname amssym.def\endcsname


\def\la{\lambda} 
\def\al{\alpha} 
\def\be{\beta}

\def\om{\omega} 
\def\({\left(} 
\def\){\right)} 
 
\def\eq{\eqalign}

\def\O#1{O\(#1\)} 
\def\abs#1{\left| #1 \right|}

\def\norm#1{\left\Vert #1 \right\Vert}

\def\klein{\eightpoint \def\smc{\small} \baselineskip=9pt}   

\def\fn#1#2{{\parindent=0.7true cm
\footnote{$^{(#1)}$}{{\klein  #2}}}}

\font\boldmas=msbm10                  
\def\Bbb#1{\hbox{\boldmas #1}}        
\def\Z{{\Bbb Z}}                        
\def\N{{\Bbb N}}                        

\def\R{{\Bbb R}}

\def\C{{\Bbb C}}
\def\A{{\Bbb A}}


\font\eightrm=cmr8                                                    
\long\def\fussnote#1#2{{\baselineskip=9pt                            
\setbox\strutbox=\hbox{\vrule height 7pt depth 2pt width 0pt}%
\eightrm                                                         
\footnote{#1}{#2}}}                                              
\font\boldmasi=msbm10 scaled 700      
\def\Bbbi#1{\hbox{\boldmasi #1}}      
\font\boldmas=msbm10                  
\def\Bbb#1{\hbox{\boldmas #1}}        
\def\Zi{{\Bbbi Z}}                      
\def\Ai{{\Bbbi A}}                      
\def\Ni{{\Bbbi N}}                      
\def\Pi{{\Bbbi P}}                      
\def\Qi{{\Bbbi Q}}                      
\def\Ri{{\Bbbi R}}



\def\dint #1 {
\quad  \setbox0=\hbox{$\disp\int\!\!\!\int$}
  \setbox1=\hbox{$\!\!\!_{#1}$}
  \vtop{\hsize=\wd1\centerline{\copy0}\copy1} \quad}

\def\drint #1 {
\qquad  \setbox0=\hbox{$\disp\int\!\!\!\int\!\!\!\int$}
  \setbox1=\hbox{$\!\!\!_{#1}$}
  \vtop{\hsize=\wd1\centerline{\copy0}\copy1}\qquad}

\def\frac#1#2{{#1\over #2}}

\def\date{\the\day.~\the\month.~\the\year}

\def\mod{\,{\rm mod}\,}  
\def\klein{\eightpoint \def\smc{\small} }

\def\at#1#2#3{{\left. \phantom{\int} #1 \right|}_{#2}^{#3}}

\def\frac#1#2{{#1\over#2}} 
\def\Int{\int\limits}  
 
\def\vol{{\rm vol}}

\hsize=16true cm     \vsize=23.2true cm

\parindent=0cm

\vbox{\vskip 2true cm}

\cen{{\bfone The lattice point discrepancy of a body of
revolution:}} \msk \cen{{\bfone Improving the lower bound by
Soundararajan's method}}\bsk \msk \cen{{\bf Manfred K\"uhleitner
and Werner Georg Nowak}}

\vbox{\vskip 1.2true cm}

\footnote{}{\klein{\it Mathematics Subject Classification } (2000): 11P21, 
11K38, 52C07.\par }

\def\SL{Soundararajan's Lemma}
\def\M{{\cal M}}
\def\MM{\widehat{\cal M}}
\def\B{{\cal B}}
\def\dB{\partial\B}
\def\La{\Lambda}
\def\dh{{\text{3\over2}}}

\def\ede{{{1\over3}}}

\def\zde{{{2\over3}}}
\def\b#1{{\bf #1}}

{\klein{\bf Abstract. } For a 
convex body $\B$ in $\Ri^3$ which is invariant under 
rotations around one coordinate axis and has a smooth boundary of 
bounded nonzero curvature, the {\it lattice point discrepancy } $P_\B(t)$ 
(number of integer points minus volume) of a linearly dilated copy 
$\sqrt{t}\B$ is estimated from below. On the basis of a recent method of 
K.~Soundararajan [16] an $\Omega$-bound is obtained that improves upon all 
earlier results of this kind.}

\vbox{\vskip 1.2true cm}

{\bf 1.~Introduction.}\quad We consider a compact convex body $\B$
in $\R^3$ which contains the origin as an inner point and assume
that its boundary $\dB$ is a $C^\infty$ surface\fn{1}{This
assumption will be made a bit more precise at the end of section
2.} with bounded nonzero Gaussian curvature throughout. For a
large real parameter $t$, we consider a linearly dilated copy
$\sqrt{t}\,\B$ of $\B$, and in particular its {\it lattice point
discrepancy}
$$ P_\B(t) := \#\(\sqrt{t}\,\B\cap\Z^3\) - \vol(\B)t^{3/2}\,.
\eqno(1.1) $$ There is a rich and very classic theory dealing with
the estimation of such quantities $P_\B(t)$, both in arbitrary
dimensions and for very special cases. An enlightening survey can
be found in E.~Kr\"atzel's monographs [8] and [9] which have to be
supplemented by M.~Huxley's book [7] where he exposed his
breakthrough in planar lattice point theory ({\it Discrete
Hardy-Littlewood method}). \ssk For our specific setting stated
above, the sharpest results read
$$ P_\B(t) = \O{t^{63/86+\ep}} \eqno(1.2) $$
and\fn{2}{For the definitions of the different $\Omega$-symbols,
cf.~Kr\"atzel [8], p.~14.}
$$ P_\B(t) = \Omega_-\(t^{1/2}(\log t)^{1/3}\)\,. \eqno(1.3) $$
These are due to W.~M\"uller [14] (who improved earlier results by
E.~Hlawka [5] and Kr\"atzel and Nowak [10], [11]), and the second
named author [15], respectively. \ssk In recent years, it has been
noted that sharper estimates are true for a body $\B$ which is
invariant under rotations around one of the coordinate axes. In
this case,
$$ P_\B(t) = \O{t^{11/16}}\,, \eqno(1.4) $$ according to
F.~Chamizo [1], and\fn{3}{By $\log_j$, $j=2,3,\dots$, we denote
throughout the $j$-fold iterated logarithm.}
$$ P_\B(t) = \Omega_-\(t^{1/2}(\log t)^{1/3}(\log_2 t)^{\ede\log2}
\exp(-c\sqrt{\log_3 t})\)\,,\quad c>0\,,\eqno(1.5) $$ as was shown
by the first named author [12], on the basis of a deep and fairly
general method of J.L.~Hafner [3]. \ssk 
Quite recently, K.~Soundararajan [16] exploited a brilliant new idea to obtain 
sharper $\Omega$-estimates in the classic circle and divisor problems. 
In the present note we will apply this ingenious new approach to 
improve\fn{4}{Note that ${\ede\log2}=0.2310\dots$ while 
$\zde(\sqrt{2}-1)=0.2761\dots $.} the lower bound of (1.5). \bsk 
{\bf Theorem.}\quad{\it Let $\B$ be a compact convex body in $\R^3$ which is 
invariant under rotations around one of the coordinate axes and contains 
$(0,0,0)$ as an inner point. Assume that its boundary $\dB$ is of class 
$C^\infty$ and has bounded nonzero Gaussian curvature throughout. Then 
$$ P_\B(t) = \Omega_-\(t^{1/2}(\log t)^{1/3}(\log_2 t)^{\zde(\sqrt{2}-1)}
(\log_3 t)^{-2/3}\)\,.$$  } \bsk\msk 
{
We remark parenthetically that still much sharper estimates are known 
for the special case that $\B$ is the unit ball $\B_0$ in $\R^3$ ({\it sphere 
problem}). In fact, Heath-Brown [4] obtained}\fn{5}{It is instructive to 
compare the numerical values of the exponents in (1.2), (1.4), and (1.6): 
${63\over86}=0.7325\dots$, ${11\over16}=0.6875$, ${21\over32}=0.65625$.} 
{
$$ P_{\B_0}(t) = \O{t^{21/32+\ep}}\,,  \eqno(1.6)$$ 
thereby improving a result of Chamizo and Iwaniec [2] and earlier classic work 
of I.M.~Vinogradov [20]. In the other direction, K.-M.~Tsang [19] showed that 
$$ P_{\B_0}(t) = \Omega_\pm\(t^{1/2}(\log t)^{1/2}\)\,, \eqno(1.7)  $$ 
the $\Omega_-$-part of this result being much older and actually due to 
G.~Szeg\"o [17].}

\bsk\msk 

\vbox{{\bf 2.~Preliminaries.} \bsk

{\bf Soundararajan's Lemma {\rm[16]}.}\quad{\it Let
$(f(n))_{n=1}^\infty$ and $(\la_n)_{n=1}^\infty$ be sequences of
non\-negative real numbers, $(\la_n)_{n=1}^\infty$ non-decreasing,
and $\sum_{n=1}^\infty f(n)<\infty$. Let $L\ge2$ be an integer and
$\La$ a positive real parameter. Suppose further that $\M$ is a
finite set of positive integers, such that $\{\la_m:\ m\in\M\ \}
\subset[\h\La,\dh\La]$. Then, for any real $T\ge2$, there exists
some $t\in[\h T, (6L)^{|\M|+1}\,T]$ with
$$ \sum_{n=1}^\infty f(n)\cos(2\pi\la_n t)\ \ge\ {1\over8}\sum_{m\in\M}f(m) -
{1\over L-1}\sum_{n:\ \la_n\le2\La}f(n) - {2\over\pi^2
T\La}\sum_{n=1}^\infty f(n)\,. $$ }} \bsk\msk

We further notice some important properties of the {\it tac
function } $H$ of a convex body $\B$ with the properties stated
above. This is defined by
$$ H(\b w) = \max_{\b x \in \B} (\b x \cdot \b w) \qquad (\b w\in\R^3)$$
where $\cdot$ denotes the standard inner product. From this the
following facts are evident: \msk

(i) $H$ is positive and homogeneous of degree 1. \ssk

(ii) There exist constants $c_2>c_1>0$, depending on $\B$, such
that for all $\b w \in \R^3$
$$  c_1\norm{\b w}\le H(\b w) \le c_2\norm{\b w}\,, \eqno(2.1) $$
where $\norm{\cdot}$ stands for the Euclidean norm throughout.
\ssk

(iii) If $\B$ is invariant with respect to rotations around the
third coordinate axis (say), then so is $H$, i.e., for all 
$(w_1,w_2,w_3)\in\R^3$, 
$$ H(w_1,w_2,w_3) = H(\sqrt{w_1^2+w_2^2},0, w_3)\,. \eqno(2.2)$$
 \msk 

{\klein It seems appropriate to say a bit more about the smoothness 
condition that $\dB$ be of class $C^\infty$. Properly speaking, this is 
supposed to mean that for every point of $\dB$ there exists a neighbourhood in 
which the corresponding portion of $\dB$ has a regular}\fn{6}{I.e., 
${\partial\b x\,\over\partial u_1}$, ${\partial\b x\,\over\partial u_2}$ are 
linearly independent.}{\klein\ parametrization $\b x=\b x(u_1,u_2)$ whose 
components are all of class $C^\infty$. However, as has been neatly worked 
out in W.~M\"uller [13], Lemmas 1 and 2, this local property implies that the 
{\it spherical map}, which sends every point of the unit sphere into that 
point of $\dB$ where the outward normal has the same direction, is globally 
one-one and $C^\infty$. 
Under these latter conditions, Hlawka's asymptotic formulas for the Fourier 
transform of the indicator function of $\B$ had been established [5], [6]. 
These in turn have been used in [15], upon which our present analysis 
will be based. \ssk For the case that $\B$ is a body of revolution (with 
respect to the $x_3$-axis, say), the conditions of our Theorem 
can be stated in a more concise form. 
 It suffices to assume that 
$$ \dB=\{\b x=(x_1,x_2,x_3)= \(\rho(\theta)\sin(\theta)\cos(\phi),\ 
\rho(\theta)\sin(\theta)\sin(\phi),\ \rho(\theta)\cos(\theta)\):\ 
0\le\theta\le\pi,\ 0\le\phi\le2\pi\ \}\,,  $$ 
where $\rho:\ \Ri\to\Ri_{>0}$ is an even function, periodic with period 
$2\pi$ and everywhere of class $C^\infty$, which satisfies throughout
$$ \rho\,\rho''-2 \rho'^2-\rho^2\ne0\,. \eqno(2.3)$$ 
\vbox{In fact, the Gaussian 
curvature $\kappa_3$ of this surface $\dB$ is readily computed as 
$$ \kappa_3(\theta)={{\d  
x_3\over\d\theta}\over\rho(\theta)\sin(\theta)}\,
{\rho(\theta)\,\rho''(\theta)-2 \rho'^2(\theta)-\rho^2(\theta)
\over(\rho^2(\theta)+\rho'^2(\theta))^2}\,. $$ We may imagine $\dB$ to be 
generated by rotation of the {\it meridian} 
$$ \{(x_1,x_3) = \(\rho(\theta)\sin(\theta),\ \rho(\theta)\cos(\theta)\):\ 
0\le\theta\le\pi\ \} $$ around the $x_3$-axis. The curvature $\kappa_2$ of the 
latter satisfies $$ \abs{\kappa_2(\theta)} = 
{\abs{\rho(\theta)\,\rho''(\theta)-2 \rho'^2(\theta)-\rho^2(\theta)}
\over(\rho^2(\theta)+\rho'^2(\theta))^{3/2}}\,. $$ 
Therefore, (2.3) guarantees the nonvanishing of $\kappa_2$, and also that of 
$\kappa_3$, since by 
geometric evidence ${\d x_3\over\d\theta}>0$ for $0<\theta<\pi$.}} 

\bsk\bsk


{\bf 3.~Proof of the Theorem.}\quad For real $t>0$, we put
$$  X=X(t) = (\log t)^{-1}\,, \ k=k(t) = t^2 \log t\,, \eqno(3.1)$$
then the {\it Borel mean-value } of the lattice rest $P_\B$ is
defined as
$$ B(t) := {1\over\Gamma(k+1)} \Int_0^\infty e^{-u} u^k
P_\B(Xu)\d u\,. \eqno(3.2) $$ We start from formula (13) in [15]:
For large $t$, and arbitrary $\ep>0$,
$$ B(t) = -{1\over2\pi}\,t\, S(t) + \O{t^{3/8+\ep}}\,, \eqno(3.3) $$
where
$$ S(t) := \sum_{0<\norm{\b m}\le t^{\ep_0}X^{-1/2}}
{\al(\b m)\over\norm{\b m}^2}\,\exp(-\h\pi^2 X H(\b m)^2)\,
\cos(2\pi H(\b m) t)\,. \eqno (3.4)$$ Here $\ep_0>0$ is a
sufficiently small constant, $\b m = (m_1,m_2,m_3)$ denotes
elements of $\Z^3$ throughout, and the coefficients $\al(\b m)$
are positive reals bounded both from above and away from 0. By
(2.2), we can rewrite this last formula as
$$ S(t) = \sum_{0<\ell+m_3^2\le t^{2\ep_0}\log t} {g(\ell,m_3)\over\ell+m_3^2}
\,\exp(-\h\pi^2 X H(\sqrt{\ell},0,m_3)^2)\,\cos(2\pi
H(\sqrt{\ell},0,m_3) t)\,, $$ with
$$ g(\ell,m_3) := \sum_{(m_1,m_2)\in\Zi^2:\atop m_1^2+m_2^2=\ell}
\al(m_1,m_2,m_3)  \asymp r(\ell)\,, \eqno(3.5) $$ $r(\ell)$ the
number of ways to write $\ell\in\N$ as a sum of two squares of
integers. \ssk In order to apply \SL, we consider a one-one map
$\b q$ of $\N_*$ onto $\N\times\Z\setminus\{(0,0)\}$, $n \mapsto
\b q(n) =(\ell,m_3)$ such that the sequence $(\la_n)_{n=1}^\infty$
defined by $$\la_n := \at{H(\sqrt{\ell},0,m_3)}{(\ell,m_3)=\b
q(n)}{}\eqno(3.6)$$ is non-decreasing\fn{7}{In other words: We
arrange the elements $(\ell,m_3)$ of
$\Ni\times\Zi\setminus\{(0,0)\}$ according to the size of the
values $H(\sqrt{\ell},0,m_3)$.}. Putting further
$$ f(n) := \at{{g(\ell,m_3)\over\ell+m_3^2}
\,\exp(-\h\pi^2 X H(\sqrt{\ell},0,m_3)^2)}{(\ell,m_3)=\b q(n)}{}
\eqno(3.7)$$ if $\ell+m_3^2\le t^{2\ep_0}\log t$, and $f(n)=0$
else, we obtain in fact
$$ S(t) = \sum_{n=1}^\infty f(n)\cos(2\pi\la_n t)\,,  $$
and are thus prepared to apply \SL. For $T\ge40$ a large real
parameter, we put $L=[(\log_2 T)^{20}]$ and assume that the set $\M$ will
be chosen such that
$$ (6L)^{|\M|+1}\le T\,. \eqno(*) $$
Then, by \SL, there exists a value $t\in[\h T,T^2]$ for which
$$ S(t) \ge\ {1\over8}\sum_{m\in\M}f(m) -
{1\over L-1}\sum_{n:\ \la_n\le2\La}f(n) - {2\over\pi^2
T\La}\sum_{n=1}^\infty f(n)\,, \eqno(3.8)$$ where $\La>0$ is a
parameter remaining to be determined. \ssk By homogeneity of the
tac-function $H$, there exist positive constants $a_2>a_1>0$ and
$a_3>a_4>0$ depending on $\B$ such that the two-dimensional
interval $[a_1, a_2]\times[a_3,a_4]$ in the $(w_1,w_3)$-plane,
say, lies between the two curves $H(w_1,0,w_3)=\h$ and
$H(w_1,0,w_3)=\dh$. Consequently, for integers $\ell>0$ and $m_3$,
the condition $(\sqrt{\ell},m_3)\in[a_1\La,
a_2\La]\times[a_3\La,a_4\La]$ always implies that
$H(\sqrt{\ell},0,m_3)\in[\h\La,\dh\La]$. \ssk Let us denote by
$\A_1$ the set of positive integers whose prime divisors are all
congruent to 1 mod 4, and by $\om(\ell)$ the number of prime divisors
of $\ell\in\N_*$. \ssk \vbox{Then we define
$$ \MM = \{(\ell,m_3)\in\N_*^2:\ a_1^2\La^2\le\ell\le a_2^2\La^2,\
a_3\La\le m_3\le a_4\La,\ \ell\in\A_1,\ \om(\ell)=[\be\log_2\La]\
\}\,,  $$ where $\be>0$ is a coefficient whose optimal choice
ultimately will be $\be=\sqrt{2}$.} \ssk Let $\M$ be the 
preimage of $\MM$ under the map $\b q$. By
construction, $\{\la_m:\ m\in\M\ \} \subset[\h\La,\dh\La]$, as
required in \SL. \ssk By (3.5) and (3.7),
$$ \eq{\sum_{m\in\M} f(m) &\gg {1\over\La^2} \sum_{a_3\La\le m_3\le
a_4\La}\ \sum_{a_1^2\La^2\le\ell\le a_2^2\La^2,\atop
\ell\in\Ai_1,\ \om(\ell)=[\be\log_2\La]} r(\ell) \cr
&\gg{1\over\La}\sum_{a_1^2\La^2\le\ell\le a_2^2\La^2,\
\ell\in\Ai_1,\ \om(\ell)=[\be\log_2\La]} r(\ell)\,, \cr}
\eqno(3.9)$$ where we have been assuming for the moment that
$$ X H(\sqrt{\ell},0,m_3)^2 \ll 1 \eqno(**) $$ for the values of
$\ell$ and $m_3$ involved. \ssk Furthermore,
$r(\ell)\ge2^{\om(\ell)}$ for $\ell\in\A_1$, and the cardinality
of $$ {\cal S}_{\La,K} := \{ \ell\in\N_*:\ a_1^2\La^2\le\ell\le
a_2^2\La^2,\ \ell\in\A_1,\ \om(\ell)=K\ \} $$ is readily estimated
after the example of Tenenbaum [18], section II.6. One may start
from the observation that, for $\Re(s)>1$, $z\in\C$ arbitrary,
$$ \sum_{n\in\Ai_1} z^{\om(n)}n^{-s} =
\prod_{p\equiv1\mod4}\(1+{z\over p^s-1}\) =
\(\zeta_{\Qi(i)}(s)\)^{z/2} G(s; z)\,, $$ where $\zeta_{\Qi(i)}$
is the Dedekind zeta-function of the Gaussian field, and $G(s; z)$
is holomorphic and bounded in every half-plane
$\Re(s)\ge\sigma_0>\h$. It follows\fn{8}{This has been noticed
already by Soundararajan [16], f.~(3.7). The authors intend to
carry out the details for the case of a general number field $\K$ 
in a forthcoming article.} that, as long as $K\ll\log_2\La$,
$$  |{\cal S}_{\La,K}| \asymp {\La^2\over\log\La}\,
{(\h\log_2\La)^{K-1}\over(K-1)!}\,. $$ With Stirling's formula in
the shape $(K-1)! \asymp K^{K-1/2}\,e^{-K}$ and the choice
$K=[\be\log_2\La]$, this gives
$$  |{\cal S}_{\La,K}| \asymp
{\La^2\over\sqrt{\log_2\La}}\,(\log\La)^{\be-1-\be\log(2\be)}\,,
$$ and thus $$ |\M|=|\MM|\asymp {\La^3\over\sqrt{\log_2\La}}\,
(\log\La)^{\be-1-\be\log(2\be)}\,,\eqno(3.10)  $$ Therefore,
recalling (3.9) and the fact that $r(\ell)\ge2^{\om(\ell)}$ for
$\ell\in\A_1$, we obtain
$$ \sum_{m\in\M}f(m) \gg {\La\over\sqrt{\log_2\La}}\,
(\log\La)^{\be-1-\be\log\be}\,.\eqno(3.11)$$ We now have to choose
$\La$ such that $(*)$ is satisfied. This is done optimally as
$$ \La = c_0 (\log T)^{1/3}(\log_2 T)^{\ede(1-\be+\be\log(2\be))}
(\log_3 T)^{-1/6}\,, \eqno(3.12)  $$ where $c_0$ is an appropriate
small constant. As a consequence, $(**)$ is verified, since
$X\ll(\log T)^{-1}$ and $H(\sqrt{\ell},0,m_3)\ll\La$ for the
values of $\ell$ and $m_3$ involved. Furthermore,
$\log\La\asymp\log_2 T$ and $\log_2\La\asymp\log_3 T$, thus
ultimately
$$  \sum_{m\in\M}f(m) \gg  (\log T)^{1/3}
(\log_2 T)^{\zde(\be-1-\be\log\be)+\ede\be\log2} (\log_3
T)^{-2/3}\,. $$ Here the second exponent is maximized for
$\be=\sqrt{2}$, and we finally obtain
$$  \sum_{m\in\M}f(m) \gg  (\log T)^{1/3}(\log_2
T)^{\zde(\sqrt{2}-1)}(\log_3 T)^{-2/3}\,.  \eqno(3.13) $$ It
remains to show that the two other terms on the right hand side of
(3.8) are small. \ssk In fact,
$$ \sum_{n:\la_n\le2\La}f(n) \ll
\sum_{0<H(\sqrt{\ell},0,m_3)\le2\La} {r(\ell)\over\ell+m_3^2}=
 \sum_{0<H(\b m)\le2\La}\norm{\b m}^{-2} \le $$
$$ \le  \sum_{0<c_1\norm{\b m}\le2\La}\norm{\b
m}^{-2} =  \sum_{1\le n\le(4/c_1^2)\La^2} {r_3(n)\over n} =
 \Int_{1-}^{(4/c_1^2)\La^2} {1\over u}
\d\(\sum_{1\le n\le u}r_3(n)\) \ \ll\La\,, $$ using integration by
parts of Stieltjes integrals and the well-known bound
$\disp\sum_{1\le n\le u}r_3(n)\ll u^{3/2}$. After division by
$L-1$, which by construction is $\asymp (\log_2 T)^{20}$, this is
small compared to the right-hand side of (3.13). \ssk
\vbox{Similarly (for the value of $t\in[\h T,T^2]$ specified by
\SL),$$ {2\over\pi^2 T\La}\sum_{n=1}^\infty f(n)\ll {1\over T\La}
\sum_{0<\norm{\b m}\le t^{\ep_0}/\sqrt{X}}\norm{\b m}^{-2} = $$
$$ = {1\over T\La} \Int_{1-}^{t^{2\ep_0}\log t}
{1\over u}\d\(\sum_{1\le n\le u}r_3(n)\)\ \ll T^{3\ep_0-1}\,. $$}
Combining the last two bounds with (3.8) and (3.3), we conclude
that for arbitrary $T\ge40$, there exists a value $t\in[\h T,T^2]$
with
$$ -B(t) \gg t (\log t)^{1/3}(\log_2 t)^{\zde(\sqrt{2}-1)}(\log_3
t)^{-2/3}\,. \eqno(3.14) $$ Let us assume that, with some
constants $C$ and $\ep_1>0$, and for all $u>0$,
$$ -P_{\B}(u) \le C + \ep_1 u^{1/2}{\cal L}(u)\,,$$ where $$ {\cal
L}(u):= (\log u)^{1/3}(\log_2 u)^{\zde(\sqrt{2}-1)} (\log_3
u)^{-2/3}\,  $$ for $u\ge20$, and ${\cal L}(u)={\cal L}(20)$ else.
By the definition (3.2) of $B(t)$, this implies that
$$ - B(t) \le C + {\ep_1\over\Gamma(k+1)} \Int_0^\infty e^{-u} u^k
(Xu)^{1/2} {\cal L}(Xu)\d u\,, $$ for all $t>0$. Estimating this
integral by Hafner's Lemma 2.3.6 in [3], we obtain
$$ -B(t)\le C + C_1 \ep_1 (kX)^{1/2} {\cal L}(kX) =  C + C_1 \ep_1
\,t\, {\cal L}(t^2)\,, $$ recalling (3.1). Together with (3.14),
this yields a positive lower bound for $\ep_1$ and thus completes
the proof of our Theorem. \qed


\vbox{\vskip 2.5true cm}   

\klein \parindent=0pt 

\cen{\bf References}  \bsk

[1] {\smc 
F.~Chamizo,} Lattice points in bodies of revolution. Acta Arith. {\bf85}, 
265-277 (1998). \ssk 

[2] {\smc F.~Chamizo \and H.~Iwaniec,} On the sphere problem. 
Rev.~Mat.~Iberoamericana {\bf 11}, 417-429 (1995). \ssk 

[3] {\smc 
J.L.~Hafner,} On the average order of a class of arithmetical functions. 
J.~Number Th. {\bf15}, 36-76 (1982).  \ssk 

[4] {\smc R.~Heath-Brown,} Lattice points in the sphere. In: Number theory in 
progress, Proc. Number Theory Conf. Zakopane 1997, eds. K.~Gy\"ory et al., 
vol.~{\bf2}, 883-892 (1999).  \ssk 

[5] {\smc 
E.~Hlawka,} \"Uber Integrale auf konvexen K\"orpern I. Monatsh.~f.~Math. 
{\bf 54}, 1-36 (1950). \ssk 
    
[6] {\smc E.~Hlawka}, {šber Integrale auf konvexen K\"orpern II}. 
Monatsh.~f.~Math. {\bf 54}, 81--99 (1950). \ssk 

[7] {\smc M.N.~Huxley}, {Area, lattice points, and exponential sums.} 
LMS Monographs, New Ser. {\bf 13}, Oxford 1996.  \ssk 

[8] {\smc E.~Kr\"atzel,} Lattice points. Berlin 1988. \ssk 

[9] {\smc E.~Kr\"atzel,} Analytische Funktionen in der Zahlentheorie. 
Stuttgart-Leipzig-Wiesbaden 2000.  \ssk 

[10] {\smc 
E.~Kr\"atzel \and W.G.~Nowak,} Lattice points in large convex bodies. 
Monatsh.~Math. {\bf112}, 61-72 (1991).  \ssk 

[11] {\smc 
E.~Kr\"atzel \and W.G.~Nowak,} Lattice points in large convex bodies II. 
Acta Arithm. {\bf62}, 232-237 (1992). \ssk 

[12] {\smc 
M.~K\"uhleitner,} Lattice points in bodies of revolution in $\Ri^3$: an 
$\Omega_-$-estimate for the error term. Arch.~Math. {\bf 74}, 234-240 (2000). 
\ssk 

[13] {\smc W.~M\"uller}, On the average order of the lattice rest of a convex 
body. Acta Arithm. {\bf 80}, 89--100 (1997).  \ssk   

[14] {\smc W.~M\"uller}, Lattice points in large convex bodies. Monatsh.~Math. 
{\bf128}, 315-330 (1999).  \ssk 

[15] {\smc W.G.~Nowak,} On the lattice rest of a convex body in $\Ri^s$, II. 
Arch.~Math. {\bf47}, 232-237 (1986). \ssk 

[16] {\smc 
K.~Soundararajan,} Omega results for the divisor and circle problems. IMRN 
{\bf36}, 1987-1998 (2003). \ssk 

[17] {\smc G.~Szeg\"o,} Beitr\"age zur Theorie der Laguerreschen Polynome, II, 
Zahlentheoretische Anwendungen. Math.~Z. {\bf25}, 388-404 (1926). \ssk 

[18] {\smc 
G.~Tenenbaum,} Introduction to analytic and probabilistic number theory. 
Cambridge 1995. \ssk 

[19] {\smc K.-M.~Tsang,} Counting lattice points in the sphere. Bull.~London 
Math.~Soc. {\bf32}, 679-688 (2000). \ssk 

[20] {\smc 
I.M.~Vinogradov,} On the number of integer points in a sphere (Russian). 
Izv.~Akad.~Nauk SSSR Ser.~Mat. {\bf27}, 957-968 (1963). \ssk 

\vbox{\vskip 1.5true cm}  

\parindent=1.5true cm

\vbox{Manfred K\"uhleitner \& Werner Georg Nowak \ssk 

Institut f\"ur Mathematik 

Department f\"ur Integrative Biologie 

Universit\"at f\"ur Bodenkultur Wien 

Peter Jordan-Stra\ss e 82 

A-1190 Wien, \"Osterreich \ssk 

E-mail: {\tt kleitner@edv1.boku.ac.at, \ nowak@mail.boku.ac.at} \ssk 

Web: http://www.boku.ac.at/math/nth.html}

\bye